\documentclass{article}

\usepackage{verbatim}
\usepackage{amsthm, amsmath}
\usepackage{indentfirst}
\usepackage{graphicx}
\usepackage{tikz}
\usepackage{stfloats}
\usepackage{enumerate}
\usepackage[margin=2.5cm]{geometry}
\usepackage[labelsep=period]{caption}
\usepackage{hyperref}

\newtheorem{lemma2026MinDegree}{Lemma}[section]
\newtheorem{theorem2026MinDegree}[lemma2026MinDegree]{Theorem}
\newtheorem{corollary2026MinDegree}[lemma2026MinDegree]{Corollary}
\newtheorem{conjecture2026MinDegree}[lemma2026MinDegree]{Conjecture}
\newtheorem{proposition2026MinDegree}[lemma2026MinDegree]{Proposition}

\begin{document}
\title{The minimum degree of minimally $k$-extendable graphs for some $k$}\date{}
\author{
 Zhiyong Gan$^{1}$\thanks{E-mail address: ganzhiyong@m.scnu.edu.cn},\quad
 Yi Huang$^{2}$\thanks{Corresponding author. E-mail address: yihuang@m.scnu.edu.cn}
\\\\1 School of Computer Science, South China Normal University, \\Guangzhou 510631, P. R. China
\\\\2 Aberdeen Institute of Data Science and Artificial Intelligence, \\South China Normal University, \\Guangzhou 510631, P. R. China
}
\maketitle
\begin{abstract}
    A graph on at least $2k+2$ vertices is $k$-extendable if every matching of size $k$ extends to a perfect matching, and it is minimally $k$-extendable if it is $k$-extendable but deleting any edge is not. Plummer proved that every $k$-extendable graph is $(k+1)$-connected. So its minimum degree is at least $k+1$. Lou and Yu conjectured that a minimally $k$-extendable graph on at most $4k$ vertices has minimum degree in $\{k+1,2k,2k+1\}$; the case $k=1$ and the claw-free case of $k=2$ have been settled previously. We determine the minimum degree of minimally $k$-extendable graphs for $k=2$ and $k=3$, thereby settling the Lou--Yu conjecture for these two values. For $k=2$ the minimum degree is in $\{3,4,5\}$. For $k=3$ it is in $\{4,7\}$ at orders 8 and 10, and in $\{4,6,7\}$ at order 12. The proofs pass to the complement graph and combine Tutte's perfect matching theorem with a classification of four-vertex and six-vertex obstructions. All predicted values are realized by explicit constructions.
\end{abstract}

\noindent\textbf{Keywords:} minimally $k$-extendable graph; perfect matching; minimum degree; matching extension.

\noindent\textbf{Mathematics Subject Classification (2020):} 05C70, 05C07.

\section{Introduction}\label{labelIntroduction}

\noindent The concept of $k$-extendable graphs was firstly introduced by Plummer \cite{articlePlummer1980} in 1980. He proved that every $k$-extendable graph is $(k+1)$-connected, in particular $\delta(G)\geq k+1$. For bipartite graphs the lower bound is sharp. Ananchuen and Caccetta \cite{articleAnanchuen1994} firstly introduced the concept of minimally $k$-extendable graphs in 1994. Lou \cite{articleLou1999} proved that every minimally $k$-extendable bipartite graph has $\delta(G)=k+1$ in 1999. In 2004, the following conjecture of Lou and Yu \cite{articleLouYu2004} predicts that for non-bipartite graphs the minimum degree can be larger, but only in a very restricted way. Gan et al. \cite{articleGanLou2019, articleGanLou2020, articleGanLouXu2020, articleGanLouZhangWen2020, articleGanXu2021, articleGanXu2022, articleGanXu2020} showed more results on $k$-extendable graphs.

\begin{conjecture2026MinDegree}[Lou and Yu \cite{articleLouYu2004}] \label{conjectureLouyu2004}
    Let $G$ be a minimally $k$-extendable graph on $\nu$ vertices with $\nu\leq 4k$. Then $\delta(G)\in\{k+1,\;2k,\;2k+1\}$.
\end{conjecture2026MinDegree}

The case $k=1$ was settled by He, Lu and Xue \cite{articleHeLuXue2025}, that is $\delta\in\{2,3\}$. Recently Guo, Lu and Zhang \cite{articleGuoLuZhang2025} settled the claw-free case of $k=2$, proving $\delta\in\{4,5\}$.

In this paper we settle the conjecture for $k=2$ and $k=3$.

\begin{theorem2026MinDegree}\label{theorem2026MinDegreeK2}
Let $G$ be a minimally 2-extendable graph on at most 8 vertices. Then $\delta(G)\in\{3,4,5\}$, and each of 3,4,5 is attained.
\end{theorem2026MinDegree}

\begin{theorem2026MinDegree}\label{theorem2026MinDegreeK3}
Let $G$ be a minimally 3-extendable graph on at most 12 vertices. Then
\begin{enumerate} [(i)]
    \item \label{theorem2026MinDegreeK31}
        if $\nu(G)=8$, then $\delta(G)\in\{4,7\}$;
    \item \label{theorem2026MinDegreeK32}
        if $\nu(G)=10$, then $\delta(G)\in\{4,7\}$;
    \item \label{theorem2026MinDegreeK33}
        if $\nu(G)=12$, then $\delta(G)\in\{4,6,7\}$.
\end{enumerate}
\end{theorem2026MinDegree}

For a 2-extendable graph, its order is even and at least 6. Hence $\nu \in \{6, 8\}$ for Conjecture \ref{conjectureLouyu2004}. Similarly, for a 3-extendable graph, we have $\nu \in \{8,10,12\}$  for Conjecture \ref{conjectureLouyu2004}. Both theorems confirm the conjecture in full. By Lemma~\ref{lemma2026MinDegreeLouYu}, $\delta=5$ at orders 10,12 is impossible. By Lemma~\ref{lemma2026MinDegree12Tight}, we have $\delta(G)\leq7$ at order 12.

The proofs share a common strategy. After deleting a suitable edge $e$, the graph $G-e$ has complement $\overline{G}\cup\{e\}$. Hence to show that every small induced subgraph still has a perfect matching, contradicting minimality. The required small-subgraph analysis is a four-vertex statement for $k=2$ (Lemma~\ref{lemma2026MinDegree8})
and a six-vertex statement for $k=3$. The constructions of Section~\ref{constructions} realize the predicted values and, in the case of $G_k$, do so for every $k$.

\section{Definitions and Notation}\label{labelDefinitions}

In this paper, we only consider finite, undirected and simple graphs. Let $G$ be a graph, $V(G)$, $E(G)$ and $\nu(G)$ denote the vertex set, the edge set and order of $G$, respectively. We also write simply $\nu$ for $\nu(G)$. Let $d_G(v)$, $\delta(G)$ and $\Delta(G)$ denote the degree, minimum degree and maximum degree of $G$, respectively. The \emph{complement} $\overline{G}$ has the same vertex set as $G$ and joins two vertices only if $G$ does not. We call two edges \emph{independent} if they are disjoint. Let $G[S]$ denote the induced subgraph for any subset $S$ of $V(G)$. Let $G\sqcup H$ denote the disjoint union and $G-H$ the graph obtained by deleting the edges of $H$ from $G$. For edge $e=uv$, let $V(e)=\{u,v\}$ and $G-e$ denote the graph with $e$ deleted. Let $K_n$ denote the complete graph on $n$ vertices, $K_{m,n}$ denote the complete bipartite graph on $m$ and $n$ vertices, $P_n$ denote the path on $n$ vertices, $C_n$ denote the cycle on $n$ vertices, and $Q_3$ denote the cube on $\{0,1\}^3$. Let $o(G)$ denote the odd components of $G$. For a graph $G$ and a vertex subset $S\subseteq V(G)$, let $G[S]$ denote the induced subgraph of $G$ with vertex set $S$. Let $\overline{G}$ denote the complement of $G$. We also write $K_{\nu(G)}-G$ for $\overline{G}$. A $t$-set is a subset of $V(G)$ of cardinality $t$. Two graphs $G$ and $H$ are \emph{isomorphic} if there exists a bijection $f: V(G) \to V(H)$ such that any two vertices $u$ and $v$ are adjacent in $G$ if and only if $f(u)$ and $f(v)$ are adjacent in $H$, and we write $G\cong H$.

A \emph{matching} in a graph is a set of pairwise disjoint edges, and it is \emph{perfect} if it covers every vertex. A \emph{$k$-matching} in a graph is a matching with exactly $k$ edges. A graph with at least $2k+2$ vertices is said to be $k$-extendable if it has a matching of size $k$, and every matching of size $k$ extends to a perfect matching. $G$ is called \emph{minimally $k$-extendable} if $G$ is $k$-extendable but $G-e$ is not $k$-extendable for any edge $e$ of $G$.

For the terminology and notation not defined in this paper, the reader is referred to \cite{bookBondy2008}.

\section{Preliminary results} \label{labelPreliminary}

In this section, we introduce some known results and show some results which will be used in the proofs of our main result.

\begin{lemma2026MinDegree}[Plummer \cite{articlePlummer1980}]\label{lemma2026MinDegreePlummer}
 Let $G$ be an $n$-extendable graph. Then the connectivity $\kappa(G)\geq n+1$.
\end{lemma2026MinDegree}

\begin{lemma2026MinDegree}[Lou and Yu \cite{articleLouYu2004}]\label{lemma2026MinDegreeLouYu}
    Let $G$ be a $k$-extendable graph on $\nu$ vertices.
\begin{enumerate}[(i)]
    \item If $k\geq\nu/4$, then $G$ is bipartite or $\kappa(G)\geq2k$. \label{lemma2026MinDegreeLouYu1}
    \item If $G$ is minimally $k$-extendable and bipartite, then $\delta(G)=k+1$. \label{lemma2026MinDegreeLouYu2}
\end{enumerate}
\end{lemma2026MinDegree}

\begin{lemma2026MinDegree}\label{lemma2026MinDegreeBasic}
    Let $G$ be a graph of even order $\nu\geq 2k+2$ and let $s=\nu-2k$.
\begin{enumerate} [(i)]
    \item $G$ is $k$-extendable if and only if for every $2k$-set $M$ such that $G[M]$ has a perfect matching, the graph $G[V(G)\setminus M]$ has a perfect matching. \label{lemma2026MinDegreeBasic1}
    \item If $G[S]$ has a perfect matching for every $s$-set $S$, then $G$ is $k$-extendable. \label{lemma2026MinDegreeBasic2}
\end{enumerate}
\end{lemma2026MinDegree}
\begin{proof}
    By the definition of $k$-extendable graphs, (i) holds. For any $k$-matching $M$, $V(G)\setminus M$ has $s$ vertices. Since $G[V(G)\setminus M]$ has a perfect matching, by the definition of $k$-extendable graphs, (ii) holds.
\end{proof}\bigskip

\begin{lemma2026MinDegree}\label{lemma2026MinDegree4}
    Let $G$ be a graph of order 4. Then $K_4-G$ contains no perfect matching if and only if $G$ contains a triangle or a copy of $K_{1,3}$.
\end{lemma2026MinDegree}
\begin{proof}
    We first prove the sufficiency. If $G$ contains a triangle, then the three vertices of the triangle are pairwise non-adjacent in $K_4-G$. If $G$ contains a $K_{1,3}$ with center $v$, then $v$ is isolated in $K_4-G$. Both cases imply that $K_4-G$ contains no perfect matching.
    
    Next, we prove the necessity. On the contrary, if $G$ contains no triangle and no $K_{1,3}$, then $\Delta(G)\leq2$. Hence $G$ is a disjoint union of paths and cycles and hence a subgraph of $C_4$. But $K_4-C_4$ is a pair of independent edges. Hence $K_4-G$ contains a perfect matching.
\end{proof}

\section{Proof of Theorem \ref{theorem2026MinDegreeK2} on minimally 2-extendable graphs}\label{sectionProof2}

\begin{lemma2026MinDegree}\label{lemma2026MinDegree8}
    Let $G$ be a graph of order 8 such that $\overline{G}$ contains no triangle and $\Delta(\overline{G})\leq 2$. Then $G$ is 2-extendable.
\end{lemma2026MinDegree}
\begin{proof}
    Since $\Delta(\overline{G})\leq 2$, $\delta(G)\geq 5$. So $G$ has two independent edges. Let $e,f$ be two independent edges and $S=V(G)\setminus(V(e)\cup V(f))$. Since $|V(G)|=8$, $S$ is a 4-set. Then $G[S]=K_4-\overline{G}[S]$, $\Delta(\overline{G}[S])\leq 2$ and $\overline{G}[S]$ is triangle-free. By Lemma~\ref{lemma2026MinDegree4}, $K_4-\overline{G}[S]$ has a perfect matching. Hence $e$ and $f$ can extend to a perfect matching. Hence $G$ is 2-extendable.
\end{proof}\bigskip

\begin{proof}[Proof of Theorem~\ref{theorem2026MinDegreeK2}]
    Let $G$ be minimally 2-extendable with $\nu(G)\leq8$. By Lemma~\ref{lemma2026MinDegreePlummer}, $\nu(G)\in\{6,8\}$ and $\delta(G)\geq3$. If $\nu(G)=6$, then $\delta(G)\leq5$ trivially. Now we consider $\nu(G)=8$.
    
    Suppose $\delta(G)\geq6$. Then $\Delta(\overline{G})\leq 1$. So $\overline{G}$ is a matching $M$ with $0\leq |M|\leq4$, and $G=K_8-M$. We delete an edge from $G$ without losing 2-extendability, contradicting minimality.

    If $|M|\leq3$, then the matching $M$ misses two vertices $u,v$. Hence $e=uv$ is an edge of $G$, and $\overline{G-e}=M\cup\{e\}$ is a matching. By Lemma~\ref{lemma2026MinDegree8}, $G-e$ is 2-extendable, contradicting minimality.

    If $|M|=4$, then $M$ is a perfect matching. Let $aa', bb'\in M$. Then $e=ab$ is an edge of $G$. In $\overline{G-e}=M\cup\{e\}$, the neighbors of $a$ are exactly $a',b$, and those of $b$ are exactly $b',a$. No vertex is adjacent to both $a$ and $b$ in $\overline{G-e}$. So $\overline{G-e}$ is triangle-free with maximum degree exactly 2. By Lemma~\ref{lemma2026MinDegree8}, $G-e$ is 2-extendable, contradicting minimality.   
\end{proof}

\begin{proposition2026MinDegree}\label{proposition2026MinDegreeK33}
    $K_{3,3}$ is minimally 2-extendable and $\delta(K_{3,3})=3$.
\end{proposition2026MinDegree}
\begin{proof}
    Let the bipartition be $A=\{a_1, a_2, a_3\}$, $B=\{b_1, b_2, b_3\}$. Easy to see $K_{3,3}$ is 2-extendable. Now for any edge of $K_{3,3}$, without loss of generality, consider $a_1b_1$. In $K_{3,3}-a_1b_1$, the independent edges $a_2b_2,a_3b_3$ leave $a_1$ and $b_1$, which are not adjacent. So $K_{3,3}-a_1b_1$ is not 2-extendable
\end{proof}\bigskip

\begin{proposition2026MinDegree}\label{proposition2026MinDegreeK6}
    $K_6$ is minimally 2-extendable and $\delta(K_6)=5$.
\end{proposition2026MinDegree}
\begin{proof}
    Let $V(K_6)=\{a_1, a_2, a_3, a_4, a_5, a_6\}$. Deleting the four endpoints of any two independent edges leaves two adjacent vertices. So $K_6$ is 2-extendable. Deleting any edge of $K_6$, without loss of generality, deleting $a_1a_2$. Then the independent edges $a_3a_4, a_5a_6$ leave $a_1, a_2$, which are not adjacent. So $K_6-a_1a_2$ is not 2-extendable.
\end{proof}\bigskip

The three values are attained by Propositions~\ref{proposition2026MinDegreeK33}, \ref{proposition2026MinDegreeK6} and Theorem~\ref{proposition2026MinDegreeMainconstr} (with $k=2$, giving $\delta=4$).

\section{Proof of Theorem \ref{theorem2026MinDegreeK3} on minimally 3-extendable graphs}\label{sectionProof3}

For $k=3$ we need six-vertex statements. We use Tutte's theorem (i.e., a graph $G$ has a perfect matching if and only if $o(G-S)\leq|S|$ for every $S\subseteq V(G)$) to prove the following lemmas.

\begin{lemma2026MinDegree}\label{lemma2026MinDegree6}
    Let $G$ be a graph of order 6 such that $\delta(G)\geq 3$. Then $G$ has a perfect matching.
\end{lemma2026MinDegree}
\begin{proof}
    Suppose not. By Tutte's theorem, there is $S\subseteq V(G)$ with $o(G-S)\geq|S|+2$. Let $C$ be an odd component of $G-S$ and $v\in C$. Since $\delta(G)\geq 3$, $3\leq d_G(v)\leq|S|+|C|-1$. So $|C|\geq 4-|S|$.
    \begin{itemize}
        \item If $|S|=0$, then every odd component has at least 5 vertices and $o(G)\geq2$ needs at least 10 vertices, a contradiction.
        \item If $|S|=1$, then every odd component has at least 3 vertices and $o(G-S)\geq 3$ needs at least 9 vertices, a contradiction.
        \item If $|S|=2$, then every odd component has at least 2 vertices. Hence at least 3 vertices. But $o(G-S)\geq 4$ needs at least 12 vertices, a contradiction.
        \item If $|S|\geq3$, then we have $o(G-S)\geq 5$, which contradicts $\nu(G)=6$.
    \end{itemize}
\end{proof}

\begin{lemma2026MinDegree}\label{lemma2026MinDegree6Sharp}
    Let $G$ be a graph of order 6, and $e=uv$ be an edge of $G$ with $\Delta(G-e)\leq 2$. Then $K_6-G$ has a perfect matching.
\end{lemma2026MinDegree}
\begin{proof}
    Note $\overline{G}=K_6-G$. Since $\Delta(G-e)\leq 2$, each vertex except $\{u,v\}$ has degree at least 3 in $\overline{G}$, and $u,v$ have degree at least 2 in $\overline{G}$. Suppose $\overline{G}$ has no perfect matching. By Tutte's theorem, there is $S\subseteq V(\overline{G})$ with $o(\overline{G}-S)\geq|S|+2$. For an odd component $C$ and $w\in C$, $d_{\overline{G}}(w)\leq|S|+|C|-1$. If $w\notin\{u,v\}$, then $|C|\geq4-|S|$, and if $w\in\{u,v\}$, then $|C|\geq3-|S|$.
    \begin{itemize}
        \item If $|S|=0$, then $o(\overline{G})\geq 2$. Since each vertex except $\{u,v\}$ has degree at least 3 in $\overline{G}$, each odd component hence has at least 5 vertices. Hence it needs 10 vertices, a contradiction.
        \item If $|S|=1$, then $o(\overline{G}-S)\geq 3$. Since each vertex except $\{u,v\}$ has degree at least 3 in $\overline{G}$, each odd component has at least 3 vertices. Hence it needs 9 vertices, a contradiction.
        \item If $|S|=2$, then $o(\overline{G}-S)\geq 4$. Since $\nu(\overline{G})=\nu(G)=6$, $o(\overline{G}-S)=4$ and each odd component is a singleton with degree at most 2 in $\overline{G}$. It contradicts the fact that each vertex except $\{u,v\}$ has degree at least 3 in $\overline{G}$.
        \item If $|S|\geq 3$, then $o(\overline{G}-S)\geq5$ while $\nu(\overline{G}-S)\leq3$, a contradiction.
    \end{itemize}
\end{proof}

\begin{corollary2026MinDegree}\label{lemma2026MinDegree6Cor}
    Let $G$ be a graph of order 6 with $\Delta(G)\leq 2$. Then $K_6-G$ has a perfect matching.
\end{corollary2026MinDegree}

We call a six-vertex graph $X$ with $\Delta(X)\leq 3$ \emph{bad} if $K_6-X$ has no perfect matching, and call $X$ \emph{good} if $K_6-X$ has a perfect matching. For a graph $G$ and a 6-set $T\subseteq V(G)$, we say that $T$ is bad or good according as $\overline{G}[T]$ is bad or good.

\begin{lemma2026MinDegree}\label{lemma2026MinDegree6BadClass}
    Let $G$ be a graph of order 6 with $\Delta(G)\leq 3$. Then $G$ is bad if and only if $G$ is isomorphic to one of $K_4\sqcup K_2$, $K_4\sqcup 2K_1$, or $K_{3,3}$.
\end{lemma2026MinDegree}
\begin{proof}
    We first prove the sufficiency. If $G=K_4\sqcup K_2$ or $G=K_4\sqcup 2K_1$, then the four vertices of $K_4$ are pairwise non-adjacent in $K_6-G$.
     So a perfect matching would have to match all four into the only two remaining vertices, which is impossible. If $G=K_{3,3}$, then $K_6-G$ is the union of two disjoint triangles and has no perfect matching.

    Now, we prove the necessity. Conversely, suppose $G$ is bad, that is, $K_6-G$ has no perfect matching. Note $\overline{G}=K_6-G$. Since $\Delta(G)\leq 3$ and $\nu(G)=6$, $\delta(\overline{G})\geq2$. By Tutte's theorem, there is $S\subseteq V(G)$ with $o(\overline{G}-S)\geq|S|+2$.
    \begin{itemize}
        \item If $|S|=0$, then $o(\overline{G})\geq 2$. Since $\delta(\overline{G})\geq2$, $G$ is exactly two disjoint triangles. So $G=K_{3,3}$.
        \item If $|S|=1$, then $o(\overline{G}-S)\geq 3$ and each odd component has at least 3 vertices. But $\nu(\overline{G}-S)=5$, a contradiction.
        \item If $|S|=2$, then $o(\overline{G}-S)\geq 4$ and $\nu(\overline{G}-S)=4$. So they are 4 singletons. Since $\delta(\overline{G})\geq2$, each singleton is adjacent to both vertices of $S$ in $\overline{G}$. So $\overline{G}$ is $K_{2,4}$ plus possibly the edge inside $S$. Taking complements, $G\cong K_4\sqcup K_2$ or $K_4\sqcup 2K_1$.
        \item If $|S|\geq3$, then $\nu(\overline{G}-S)\leq 3$, contradicting $o(\overline{G}-S)\geq 5$.
    \end{itemize}
\end{proof}

\subsection{The case of order 8 (Part (\ref{theorem2026MinDegreeK31}) of Theorem~\ref{theorem2026MinDegreeK3})}

Let $G$ be a 3-extendable graph of order 8. We show $\delta(G)\in\{4,7\}$.
\begin{proof}[Proof of Part \eqref{theorem2026MinDegreeK31} of Theorem~\ref{theorem2026MinDegreeK3}]
     Note $\delta(G)=7-\Delta(\overline{G})$. If $\overline{G}$ has no edges, then $G\cong K_8$ and $\delta(G)=7$. Now assume that $\overline{G}$ has an edge $uv$, i.e., $u,v$ are disjoint in $G$. Let $H=G-\{u,v\}$. Since $G$ is 3-extendable, by Lemma~\ref{lemma2026MinDegreeBasic}~\eqref{lemma2026MinDegreeBasic1},  $H=K_6-\overline{H}$ has no perfect matching. By Corollary~\ref{lemma2026MinDegree6Cor}, we have $\Delta(\overline{G})\geq 3$. Hence $\delta(G)=7-\Delta(\overline{G})\leq 4$. By Lemma~\ref{lemma2026MinDegreePlummer}, $\delta(G)\geq4$. So $\delta(G)=4$. So the part \eqref{theorem2026MinDegreeK31} of Theorem~\ref{theorem2026MinDegreeK3} holds.
\end{proof}\bigskip
\noindent\textbf{Remark.}\emph{ Part~\eqref{theorem2026MinDegreeK31} in fact holds for every 3-extendable graph of order 8, not only the minimally 3-extendable ones, since the proof uses only 3-extendability.}

\begin{proposition2026MinDegree}\label{proposition2026MinDegreeK44}
    $K_{4,4}$ is minimally 3-extendable and $\delta(K_{4,4})=4$.
\end{proposition2026MinDegree}
\begin{proof}
    Let the bipartition be $A=\{a_1, a_2, a_3, a_4\}$, $B=\{b_1, b_2, b_3, b_4\}$. A 3-matching of $K_{4,4}$ uses three vertices from each part and leaves one vertex from each, which are adjacent. Hence $K_{4,4}$ is 3-extendable.  Deleting any edge of $K_{4,4}$, without loss of generality, deleting $a_1b_1$, the 3-matching $a_2b_2,a_3b_3,a_4b_4$ leaves $a_1,b_1$ disjoint. So $K_{4,4}-a_1b_1$ is not 3-extendable.
\end{proof}\bigskip

\begin{proposition2026MinDegree}\label{proposition2026MinDegreeK8}
    $K_8$ is minimally 3-extendable and $\delta(K_8)=7$.
\end{proposition2026MinDegree}
\begin{proof}
    Let $V(K_8)=\{a_1, a_2, \ldots, a_8\}$. Deleting the six endpoints of any two independent edges leaves two adjacent vertices. So $K_8$ is 3-extendable. Deleting any edge of $K_8$, without loss of generality, deleting $a_1a_2$. Then the independent edges $a_3a_4, a_5a_6, a_7a_8$ leave $a_1, a_2$, which are not adjacent. So $K_8-a_1a_2$ is not 3-extendable.

\end{proof}

\subsection{The case of order 10 (Part (\ref{theorem2026MinDegreeK32}) of Theorem~\ref{theorem2026MinDegreeK3})}

Let $G$ be a minimally 3-extendable graph on 10 vertices. We show $\delta(G)\in\{4,7\}$.
\begin{proof}[Proof of Part \eqref{theorem2026MinDegreeK32} of Theorem~\ref{theorem2026MinDegreeK3}]
    Suppose $\delta(G)\geq8$. Then $\Delta(\overline{G})=9-\delta(G)\leq 1$. So $\overline{G}$ is a matching. For any edge $e$ of $G$, we have $\overline{G-e}=\overline{G}\cup\{e\}$. Hence $\Delta(\overline{G-e})\leq2$. For any 4-set $S$, the graph $\overline{G-e}[S]$ is a matching plus at most one edge. So $|E(\overline{G-e}[S])|\leq3$ and $\Delta(\overline{G-e}[S])\leq2$. Clearly, $\overline{G-e}[S]$ is neither a triangle nor a $K_{1,3}$. By Lemma~\ref{lemma2026MinDegree4}, $(G-e)[S]=K_4-\overline{G-e}[S]$ has a perfect matching. By Lemma~\ref{lemma2026MinDegreeBasic}~\eqref{lemma2026MinDegreeBasic2}, $G-e$ is 3-extendable, contradicting minimality. Hence $\delta(G)\leq7$.
    
    By Lemma~\ref{lemma2026MinDegreePlummer}, $\delta(G)\geq4$. By Lemma~\ref{lemma2026MinDegreeLouYu}~\eqref{lemma2026MinDegreeLouYu1}, $G$ is bipartite or $\kappa(G)\geq 6$. If $G$ is bipartite, by Lemma~\ref{lemma2026MinDegreeLouYu}~\eqref{lemma2026MinDegreeLouYu2}, $\delta(G)=4$. If $\delta(G)=5$, then $\kappa(G)\leq5<6=2k$. By Lemma~\ref{lemma2026MinDegreeLouYu}~\eqref{lemma2026MinDegreeLouYu1}, $G$ is bipartite, and Lemma~\ref{lemma2026MinDegreeLouYu}~\eqref{lemma2026MinDegreeLouYu2} gives $\delta(G)=4$, a contradiction. Now assume $\delta(G)=6$. So $\Delta(\overline{G})=9-\delta(G)=3$. Choose a vertex $v$ with $d_{\overline{G}}(v)=3$ and let $N_{\overline{G}}(v)=\{a,b,c\}$. Hence the 4-set $S=\{v,a,b,c\}$ contains a $K_{1,3}$ in $\overline{G}$. By Lemma~\ref{lemma2026MinDegree4}, the graph $G[S]=K_4-\overline{G}[S]$ has no perfect matching. Since $G$ is 3-extendable, by Lemma~\ref{lemma2026MinDegreeBasic}~\eqref{lemma2026MinDegreeBasic1}, $G[W]$ has no perfect matching, where $W=V(G)\setminus S$ has six vertices. Note $G[W]=K_6-\overline{G}[W]$. Since $\Delta(\overline{G}[W])\leq3$, by Lemma~\ref{lemma2026MinDegree6BadClass}, $\overline{G}[W]\cong K_4\sqcup 2K_1$, $K_4\sqcup K_2$ or $K_{3,3}$.

    \noindent\textbf{Case 1: $\overline{G}[W]\cong K_{3,3}$.}

    Since $\Delta(\overline{G})=3$, all the neighbors of the six vertices of $W$ in $\overline{G}$ are inside $W$, i.e., no vertex of $S$ is adjacent to any vertex of $W$ in $\overline{G}$. So for each of $a,b,c$, all of its neighbors in $\overline{G}$ are inside $S$. We split according to $\overline{G}[\{a,b,c\}]$.

    \smallskip
    \noindent\textit{Subcase 1.1: $\overline{G}[\{a,b,c\}]$ is a triangle.}
    
    Choose $x\in W$ and let $B=\{a,b,c,x\}$. The three vertices $a,b,c$ form a triangle in $\overline{G}$. By Lemma~\ref{lemma2026MinDegree4}, the graph $G[B]=K_4-\overline{G}[B]$ has no perfect matching. Since $V(G)\setminus B$ is a 6-set, $\overline{G}[V(G)\setminus B]=K_1\sqcup K_{2,3}$. The vertex $v$ is isolated in $\overline{G}$. The set $W\setminus\{x\}$ is $K_{3,3}$ with one vertex removed, i.e., $K_{2,3}$. Therefore $G[V(G)\setminus B]=K_6-(K_1\sqcup K_{2,3})$ has a perfect matching. For example, match $v$ with a vertex of the three-vertex part of the $K_{2,3}$, match the remaining two vertices of that part, and match the two vertices of the two-vertex part. So a 6-set has a perfect matching, but its complementary 4-set does not. This contradicts Lemma~\ref{lemma2026MinDegreeBasic}~\eqref{lemma2026MinDegreeBasic1}.

    \smallskip
    \noindent\textit{Subcase 1.2: $\overline{G}[\{a,b,c\}]$ is not a triangle.}
    
    Choose $p\in W$ and let $B$ be the 4-set formed by $p$ together with its three neighbors in $W$. Then $\overline{G}[B]\cong K_{1,3}$. By Lemma~\ref{lemma2026MinDegree4}, $G[B]$ has no perfect matching. The complement $V(G)\setminus B=S\cup(W\setminus B)$ is a 6-set. Since $a,b,c$ are adjacent to $v$ in $\overline{G}$, $v$ is isolated within $S$ in $G$. So in any perfect matching of $G[V(G)\setminus B]$, the vertex $v$ must be matched to a vertex of $W\setminus B$. Both vertices of $W\setminus B$ are adjacent to each vertex of $S$ in $G$. Match $v$ with a vertex of $W\setminus B$. The remaining vertex of $W\setminus B$ (say $u$) and $a,b,c$ is a 4-set. $u$ is adjacent to all of $a,b,c$ in $G$. Since $\overline{G}[\{a,b,c\}]$ is not a triangle, $G[\{a,b,c\}]$ contains at least one edge. Match that edge, and the leftover vertex of $\{a,b,c\}$ with $u$. Then it is a perfect matching of $G[V(G)\setminus B]$. Again a 6-set has a perfect matching but its complementary 4-set does not. This contradicts Lemma~\ref{lemma2026MinDegreeBasic}~\eqref{lemma2026MinDegreeBasic1}.

    Both subcases contradict 3-extendability.

    \noindent\textbf{Case 2: $\overline{G}[W]\cong K_4\sqcup 2K_1$ or $K_4\sqcup K_2$.}
    
    Let $P\subseteq W$ be the four vertices spanning the $K_4$, and let $P=\{p_1,p_2,p_3,p_4\}$. Since every vertex of $P$ in $\overline{G}$ has degree 3 inside $P$, no vertex of $P$ is adjacent to a vertex outside $P$ in $\overline{G}$. Hence $P$ is an independent set of $G$. Let $W\setminus P=\{r,s\}$.

    We claim that in a 3-extendable graph on 10 vertices, a 4-set $T$ has a perfect matching in $G$ if and only if the complementary 6-set $V(G)\setminus T$ has one. By Lemma~\ref{lemma2026MinDegreeBasic}~\eqref{lemma2026MinDegreeBasic1}, the sufficiency holds. Now we prove the necessity. Suppose $G[T]$ has a perfect matching. Then $G[V(G)\setminus T]$ has at least an edge $f$, otherwise $G$ would have no perfect matching, contradicting 3-extendability. The two matching edges of $G[T]$ and $f$ form a 3-matching. It extends to a perfect matching of $G$, and there are three edges of this perfect matching lie inside $V(G)\setminus T$.

    Now consider the 6-set $\{v,s\}\cup P$. The four vertices of $P$ are independent in $G$ and can be matched only through $v$ and $s$, which cover at most two of them. So $G[\{v,s\}\cup P]$ has no perfect matching. By the claim, $G[\{a,b,c,r\}]$ has no perfect matching. By Lemma~\ref{lemma2026MinDegree4}, $\overline{G}[\{a,b,c,r\}]$ contains a triangle or a $K_{1,3}$. The same holds for $\{a,b,c,s\}$.

    Since $\Delta(\overline{G})=3$, each of $a,b,c$ has degree at most 3 in $\overline{G}$, and each is adjacent to $v$ in $\overline{G}$ (because $N_{\overline{G}}(v)=\{a,b,c\}$) and not to any vertex of $P$. So $a$ has at most two neighbors in $\overline{G}$ among $b,c,r,s$; likewise $b$ among $a,c,r,s$ and $c$ among $a,b,r,s$. Also $N_{\overline{G}}(r)\subseteq\{a,b,c,s\}$ and $N_{\overline{G}}(s)\subseteq\{a,b,c,r\}$. So, a $K_{1,3}$ inside $\{a,b,c,r\}$ can only be centered at $r$ with $N_{\overline{G}}(r)=\{a,b,c\}$. Similarly for $\{a,b,c,s\}$.

    \smallskip
    \noindent\textit{Subcase 2.1: $r$ and $s$ are adjacent in $\overline{G}$.}

    Then $r$ has at most two neighbors among $a,b,c$. So $\overline{G}[\{a,b,c,r\}]$ contains a triangle and no $K_{1,3}$. Similarly for $\overline{G}[\{a,b,c,s\}]$. Suppose the triangle contains $r$, say $\{a,b,r\}$. Then $a$ and $b$ already have the neighbors $v,r$ and each other. So neither is adjacent to $c$ or $s$. Then $\{a,b,c,s\}$ contains no triangle, a contradiction. Hence the triangle is $\{a,b,c\}$. Then $a,b,c$ have no neighbors in $\{r,s\}$, and $\overline{G}\cong K_4\sqcup K_4\sqcup K_2$.

    \smallskip
    \noindent\textit{Subcase 2.2: $r$ and $s$ are not adjacent in $\overline{G}$.}

    If the triangle for $\{a,b,c,r\}$ contains $r$, then by the same degree count, $\{a,b,c,s\}$ contains no triangle and no $K_{1,3}$, a contradiction. So the triangle is $\{a,b,c\}$ or the $K_{1,3}$ is centered at $r$; likewise for $\{a,b,c,s\}$. If the triangle is $\{a,b,c\}$, then $a,b,c$ have no neighbors in $\{r,s\}$, and $\overline{G}\cong K_4\sqcup K_4\sqcup 2K_1$. Otherwise $N_{\overline{G}}(r)=N_{\overline{G}}(s)=\{a,b,c\}$. Then $a,b,c$ have neighbors exactly $v,r,s$. Hence $\overline{G}\cong K_4\sqcup K_{3,3}$.

    \smallskip
    If $\overline{G}\cong K_4\sqcup K_4\sqcup K_2$ or $\overline{G}\cong K_4\sqcup K_4\sqcup 2K_1$, then $G$ is the join of the independent sets $\{v,a,b,c\}$ and $P$, together with $r,s$, which are adjacent in $G$ only in the second graph. Choose $\{ra,sb,cp_1\}$ to be a 3-matching. In an extension, $v$ must take one of $p_2,p_3,p_4$. But the remaining two vertices of $P$ are not adjacent. If $\overline{G}\cong K_4\sqcup K_{3,3}$, then $\{v,r,s\}$ and $\{a,b,c\}$ are two triangles of $G$. Choose $\{rs,ab,cp_1\}$ to be a 3-matching. In an extension, $v$ must take one of $p_2,p_3,p_4$. But the remaining two vertices of $P$ are not adjacent.
    
    Both subcases contradict 3-extendability.

    So the part \eqref{theorem2026MinDegreeK32} of Theorem~\ref{theorem2026MinDegreeK3} holds.
\end{proof}

\subsection{The case of order 12 (Part (\ref{theorem2026MinDegreeK33}) of Theorem~\ref{theorem2026MinDegreeK3})}

We first prove an upper bound $\delta(G)\leq8$.

\begin{lemma2026MinDegree}\label{lemma2026MinDegree12Upper}
    If $G$ is a minimally 3-extendable graph on 12 vertices, then $\delta(G)\leq8$.
\end{lemma2026MinDegree}
\begin{proof}
    Suppose $\delta(G)\geq9$. Then $\Delta(\overline{G})\leq2$. For any edge $e$ of $G$, we have $\overline{G-e}=\overline{G}\cup\{e\}$. Let $S$ be any 6-set. If $e$ does not have both endpoints in $S$, then $\overline{G-e}[S]=\overline{G}[S]$ and $\Delta(\overline{G-e}[S])\leq2$. By Corollary~\ref{lemma2026MinDegree6Cor}, $K_6-\overline{G-e}[S]$ has a perfect matching. If both endpoints of $e$ lie in $S$, then $\overline{G-e}[S]-e=\overline{G}[S]$ and $\Delta(\overline{G}[S])\leq2$. By Lemma~\ref{lemma2026MinDegree6Sharp}, $K_6-\overline{G-e}[S]$ has a perfect matching. By Lemma~\ref{lemma2026MinDegreeBasic}~\eqref{lemma2026MinDegreeBasic2}, $G-e$ is 3-extendable. This contradicts minimality. So $\delta(G)\leq8$.
\end{proof}\bigskip

Now we improve $\delta\leq8$ to $\delta\leq7$. Recall that a \emph{branch} of $\overline{G}$ is a set of four or six vertices on which $\overline{G}$ induces a $K_4$ or a $K_{3,3}$. A vertex lies in a branch only if all three of its neighbors in $\overline{G}$ lie in that branch.

\begin{lemma2026MinDegree}\label{lemma2026MinDegreeStruct12}
    Let $G$ be a 3-extendable graph on 12 vertices with $\Delta(\overline{G})\leq3$. Then either $\overline{G}\cong K_{3,3}\sqcup K_{3,3}$, or $\overline{G}$ has no $K_4$-branch and no $K_{3,3}$-branch.
\end{lemma2026MinDegree}
\begin{proof}
    Let $A=\{T\subseteq V(G):|T|=6,\;\overline{G}[T]\text{ is bad}\}$. By Lemma~\ref{lemma2026MinDegreeBasic}~\eqref{lemma2026MinDegreeBasic1}, $G$ is 3-extendable if and only if $T\in A$ implies $V(G)\setminus T\in A$. If $T\in A$, then $\overline{G}[T]$ is bad with $\Delta(\overline{G}[T])\leq\Delta(\overline{G})\leq3$. By Lemma~\ref{lemma2026MinDegree6BadClass}, $\overline{G}[T]\cong K_4\sqcup2K_1$, $K_4\sqcup K_2$ or $K_{3,3}$. In the first two cases, since $\Delta(\overline{G})\leq3$, the four vertices of the $K_4$ have degree 3 in $\overline{G}$. So they form a $K_4$-branch. In the third case, $T$ is a $K_{3,3}$-branch. Thus every $T\in A$ is either itself a $K_{3,3}$-branch of $\overline{G}$ or contains the four vertices of a $K_4$-branch of $\overline{G}$.

    We claim that distinct branches are vertex-disjoint. Suppose two branches share a vertex, and consider the possible types.
    \begin{itemize}
        \item First, a $K_4$-branch and a $K_{3,3}$-branch cannot share a vertex. Indeed, if $v$ lies in a $K_4$-branch, then all three neighbors of $v$ in $\overline{G}$ lie in that branch and are pairwise adjacent in $\overline{G}$ (they are the other three vertices of the $K_4$). Hence $N_{\overline{G}}(v)$ is a clique with three vertices. If $v$ also lies in a $K_{3,3}$-branch, then all three neighbors of $v$ in $\overline{G}$ lie in that branch and form the part opposite to $v$. Hence $N_{\overline{G}}(v)$ is an independent set with three vertices. Thus $N_{\overline{G}}(v)$ would be simultaneously a clique and an independent set, a contradiction.

        \item Second, two distinct $K_4$-branches cannot share a vertex. If $v$ lies in a $K_4$-branch $B$, then the other three vertices of $B$ are exactly the three neighbors of $v$ in $\overline{G}$ since $\Delta(\overline{G})\leq3$. Hence $B=\{v\}\cup N_{\overline{G}}(v)$, which is uniquely determined by $v$. Any two $K_4$-branches containing $v$ must therefore coincide, a contradiction.

        \item Third, two distinct $K_{3,3}$-branches cannot share a vertex. Let $C$ and $C'$ be two $K_{3,3}$-branches. Suppose $v$ lies in both $C$ and $C'$. Since $\Delta(\overline{G})\leq3$, the three neighbors of $v$ in $\overline{G}$ all lie in the part opposite to $v$ in $C$. Therefore $N_{\overline{G}}(v)$ is exactly that part. The same holds in $C'$. So $N_{\overline{G}}(v)$ is the same 3-set in both branches, and the two branches share not only $v$ but also all of $N_{\overline{G}}(v)$. Choose a common neighbor $w\in N_{\overline{G}}(v)$. In each of $C$ and $C'$, the vertex $w$ is adjacent in $\overline{G}$ to the three vertices of the part containing $v$. These two parts both contain $v$, and they cannot be equal (otherwise $C=C'$). Their union therefore has size at least 4. So $w$ has at least four neighbors in $\overline{G}$, contradicting $\Delta(\overline{G})\leq3$.
    \end{itemize}
    
    Thus distinct branches are vertex-disjoint.

    If $\overline{G}$ has neither a $K_4$-branch nor a $K_{3,3}$-branch, then the lemma holds. We may therefore assume that $\overline{G}$ has a $K_4$-branch or a $K_{3,3}$-branch.
    \begin{itemize}
        \item If $\overline{G}$ has a $K_4$-branch $C$, then since distinct branches are vertex-disjoint, any other branch is contained in the eight vertices outside $C$ and each contains four or six vertices. So there are at most two of them. Choose $x,y$ outside $C$ such that each branch outside $C$ contains at least one of $x,y$. Let $T=C\cup\{x,y\}$. So $V(G)\setminus T$ contains no whole branch, i.e., $V(G)\setminus T\notin A$.. Since the four vertices of $C$ are pairwise non-adjacent in $K_6-\overline{G}[T]$ with only $x,y$ to match them, $T\in A$. This breaks closure.
        \item If $\overline{G}$ has a $K_{3,3}$-branch, then it has exactly two or exactly one. In the first case, $\overline{G}\cong K_{3,3}\sqcup K_{3,3}$. In the latter case, let $T$ be that branch. Then $T\in A$ but $V(G)\setminus T$ contains no whole branch. This again breaks closure.
    \end{itemize}
\end{proof}

\begin{lemma2026MinDegree}\label{lemma2026MinDegreeMarkedge}
    Let $G$ be a 3-extendable graph on 12 vertices with $\Delta(\overline{G})\leq3$. Suppose $\overline{G}$ has neither a $K_4$-branch nor a $K_{3,3}$-branch. Let $e=uv$ be an edge of $G$ and $S$ a 6-set containing $u,v$ such that $\overline{G-e}[S]$ is bad and $\overline{G-e}[V(G)\setminus S]$ is good. Then $\overline{G}[S]$ is good, $d_{\overline{G}[S]}(u)\geq2$, $d_{\overline{G}[S]}(v)\geq2$, and at least one of the following holds:
    \begin{enumerate}[(i)]
        \item $\overline{G}[\{u,v,p,q\}]\cong K_4-uv$ for some $p,q\in S$;\label{lemma2026MinDegreeMarkedge1}
        \item $N_{\overline{G}[S]}(u)=N_{\overline{G}[S]}(v)$;\label{lemma2026MinDegreeMarkedge2}
        \item $\overline{G}[S]\cong K_{3,3}-uv$.\label{lemma2026MinDegreeMarkedge3}
    \end{enumerate}
\end{lemma2026MinDegree}
\begin{proof}
    If $\overline{G}[S]$ is bad, then $G[S]=K_6-\overline{G}[S]$ has no perfect matching. By Lemma~\ref{lemma2026MinDegreeBasic}~\eqref{lemma2026MinDegreeBasic1}, $G[V(G)\setminus S]$ also has none. Then $\overline{G}[V(G)\setminus S]$ is bad. But $\overline{G-e}[V(G)\setminus S]\supseteq\overline{G}[V(G)\setminus S]$. Hence $\overline{G-e}[V(G)\setminus S]$ is bad as well, contradicting that it is good. So $\overline{G}[S]$ is good.

    Since $\Delta(\overline{G}[S])\leq3$ and $\overline{G-e}[S]=\overline{G}[S]\cup\{uv\}$, only $u$ or $v$ can have degree 4 in $\overline{G-e}[S]$. So $\Delta(\overline{G-e}[S])\leq4$.

    If $\Delta(\overline{G-e}[S])\leq3$, since $\overline{G-e}[S]$ is bad, then by Lemma~\ref{lemma2026MinDegree6BadClass}, $\overline{G-e}[S]$ is $K_4\sqcup2K_1$, $K_4\sqcup K_2$ or $K_{3,3}$. In the first two cases, $\overline{G-e}[S]$ contains a $K_4$ and the edge $uv$ must lie in that $K_4$. Otherwise, $\overline{G}[S]$ would still contain the $K_4$, contradicting the absence of a $K_4$-branch. Thus (i) holds and $d_{\overline{G}[S]}(u)=2$. In the last case, $\overline{G}[S]=\overline{G-e}[S]-uv\cong K_{3,3}-uv$, (iii) holds.

    Now assume $\Delta(\overline{G-e}[S])=4$. Since only $u$ or $v$ have degree 4 in $\overline{G-e}[S]$, without loss of generality, let $u$ have degree 4 in $\overline{G-e}[S]$. Then $u$ has degree 1 in $K_6-\overline{G-e}[S]$, let $w$ be the neighbor of $u$. So $K_6-\overline{G-e}[S]$ has a perfect matching if and only if $K_4-\overline{G-e}[S\setminus\{u,w\}]$ does. Since $\overline{G-e}[S]$ is bad, $K_4-\overline{G-e}[S\setminus\{u,w\}]$ has no perfect matching. By Lemma~\ref{lemma2026MinDegree4}, $\overline{G-e}[S\setminus\{u,w\}]$ contains a triangle or a $K_{1,3}$. Note that each vertex of $S\setminus\{u,w\}$ is the neighbor of $u$ in $\overline{G-e}[S\setminus\{u,w\}]$. A triangle not containing $v$ would, together with $u$, span a $K_4$ in $\overline{G}[S]$. This contradicts the absence of a $K_4$-branch. So the triangle is $\{v,p,q\}$ and $\overline{G}[\{u,v,p,q\}]\cong K_4-uv$, (i) holds and $d_{\overline{G}[S]}(u)=3$. A $K_{1,3}$ in $\overline{G-e}[S\setminus\{u,w\}]$ must have center $v$. Otherwise its center $c\neq v$ would have $d_{\overline{G}[S]}(c)\geq4$. Let $p,q,r$ be its leaves. Both $u$ and $v$ are adjacent in $\overline{G}[S]$ to $p,q,r$. So $N_{\overline{G}[S]}(u)=N_{\overline{G}[S]}(v)=\{p,q,r\}$, (ii) holds and $d_{\overline{G}[S]}(u)=3$.
\end{proof}\bigskip

\begin{corollary2026MinDegree}\label{lemma2026MinDegreeNodeg2}
    Let $G$ be a 3-extendable graph on 12 vertices with $\Delta(\overline{G})\leq3$, and suppose $\overline{G}$ has neither a $K_4$-branch nor a $K_{3,3}$-branch. Let $e=uv$ be an edge of $G$. If $G-e$ is not 3-extendable, then $d_{\overline{G}}(u)\geq2$ and $d_{\overline{G}}(v)\geq2$.
\end{corollary2026MinDegree}
\begin{proof}
    By Lemma~\ref{lemma2026MinDegreeBasic}~\eqref{lemma2026MinDegreeBasic1}, there is a 6-set $T$ with $(G-e)[T]$ having a perfect matching but $(G-e)[V(G)\setminus T]$ not. Let $S=V(G)\setminus T$. Then $\overline{G-e}[S]$ is bad and $\overline{G-e}[T]$ is good. If one of $u,v$ is outside $S$, $e$ does not lie in $\overline{G-e}[S]$. Then $\overline{G}[S]$ is bad, i.e., $G[S]=K_6-\overline{G}[S]$ has no perfect matching. By Lemma~\ref{lemma2026MinDegreeBasic}~\eqref{lemma2026MinDegreeBasic1}, $G[V(G)\setminus S]=G[T]$ has no perfect matching. Hence $\overline{G}[T]$ is bad. Since $\overline{G-e}[T]\supseteq\overline{G}[T]$, $\overline{G-e}[T]$ is bad, a contradiction. Thus $u,v\in S$. By Lemma~\ref{lemma2026MinDegreeMarkedge}, $d_{\overline{G}}(u)\geq d_{\overline{G}[S]}(u)\geq2$ holds and the same for $v$.
\end{proof}\bigskip

\begin{lemma2026MinDegree}\label{lemma2026MinDegree12Tight}
    If $G$ is a minimally 3-extendable graph on 12 vertices, then $\delta(G)\leq7$.
\end{lemma2026MinDegree}
\begin{proof}
    Suppose $\delta(G)\geq 8$. Then $\Delta(\overline{G})\leq 3$. By Lemma~\ref{lemma2026MinDegreeStruct12}, either $\overline{G}\cong K_{3,3}\sqcup K_{3,3}$ or $\overline{G}$ has no $K_4$-branch and no $K_{3,3}$-branch.

    \noindent\textbf{Case 1: $\overline{G}\cong K_{3,3}\sqcup K_{3,3}$.}
    
    Then $G$ is the join of two copies of $2K_3$. Let $X,Y$ denote two copies of $2K_3$, respectively. Each is the disjoint union of two triangles. So $G=X\bigvee Y$. By Lemma~\ref{lemma2026MinDegree6BadClass}, the only bad 6-sets of $\overline{G}$ are $\overline{G}[X]$ and $\overline{G}[Y]$. Choose an edge $e=xy$ with $x\in X,y\in Y$. By Lemma~\ref{lemma2026MinDegreeBasic}~\eqref{lemma2026MinDegreeBasic1}, $G-e$ is 3-extendable if and only if for every 6-set $T$ such that $(G-e)[T]$ has no perfect matching, the graph $(G-e)[V(G)\setminus T]$ has no perfect matching either. We have $\overline{G-e}[X]=\overline{G}[X]\cong K_{3,3}$ and $\overline{G-e}[Y]=\overline{G}[Y]\cong K_{3,3}$. So neither $(G-e)[X]$ nor $(G-e)[Y]$ has a perfect matching, and $Y=V(G)\setminus X$.
    
    For any other 6-set $T$, if $T$ does not contain both $x$ and $y$, then $\overline{G-e}[T]=\overline{G}[T]$ is good, that is, $(G-e)[T]$ has a perfect matching. If $T$ contains both $x$ and $y$, write the four triangle-parts as $X_a,X_b\subseteq X$ and $Y_a,Y_b\subseteq Y$, with $x\in X_a$ and $y\in Y_a$. Let $a=|T\cap X_a|$, $b=|T\cap X_b|$, $c=|T\cap Y_a|$, $d=|T\cap Y_b|$. Then $a+b+c+d=6$ with $a,c\geq1$. The graph $(G-e)[T]$ is the join of the cliques $K_a\sqcup K_b$ and $K_c\sqcup K_d$. This join always has a perfect matching. Let $s_X=a+b$ and $s_Y=c+d$. We match $m$ vertices across $K_a\sqcup K_b$ and $K_c\sqcup K_d$, the rest inside their cliques. We need $m\equiv s_X\pmod2$ and $0\leq m\leq\min(s_X,s_Y)$. If $s_X$ is even and all four cliques are even, then take $m=0$. If $s_X$ is even but some clique is odd, each side has an even number of odd cliques, then take $m=2$. If $s_X$ is odd, each side has exactly one odd clique, then take $m=1$. The required cross edges are present in the join. Thus the only 6-sets $T$ for which $(G-e)[T]$ has no perfect matching are $X$ and $Y$. If $T=X$ then $V(G)\setminus T=Y$, and if $T=Y$ then $V(G)\setminus T=X$. In both cases the complement has no perfect matching either. By Lemma~\ref{lemma2026MinDegreeBasic}~\eqref{lemma2026MinDegreeBasic1}, $G-e$ is 3-extendable. Hence $G$ is not minimally 3-extendable, a contradiction.

    \noindent\textbf{Case 2: $\overline{G}$ has no $K_4$-branch and no $K_{3,3}$-branch.}

    We find an edge $e$ of $G$ with $G-e$ still 3-extendable.
    
    \begin{itemize}
        \item If $\delta(\overline{G})\leq1$, then choose $w$ with $d_{\overline{G}}(w)\leq1$ and a vertex $x$ not adjacent to $w$ in $\overline{G}$. Then $e=wx$ is an edge of $G$. By Corollary~\ref{lemma2026MinDegreeNodeg2}, $G-e$ is 3-extendable. Otherwise we have $d_{\overline{G}}(w)\geq2$, a contradiction.

        \item If $\delta(\overline{G})=2$, then choose $w$ with $N_{\overline{G}}(w)=\{p,q\}$. Since $\Delta(\overline{G})\leq3$, each of $p,q$ has at most two neighbors other than $w$ in $\overline{G}$. Hence in $\overline{G}$ at most four of the nine vertices not adjacent to $w$ are adjacent to $p$ or $q$, and at least five of them are not adjacent to $p$ or $q$. Choose such an $x$ and let $e=wx$ in $G$. Then $N_{\overline{G}}(w)\cap N_{\overline{G}}(x)=\emptyset$. Suppose $G-e$ is not 3-extendable. By Lemma~\ref{lemma2026MinDegreeBasic}~\eqref{lemma2026MinDegreeBasic1}, there is a 6-set $T$ such that $(G-e)[T]$ has a perfect matching but $(G-e)[V(G)\setminus T]$ does not. Let $S=V(G)\setminus T$. Then $\overline{G-e}[S]$ is bad and $\overline{G-e}[T]$ is good. We claim that $w,x\in S$. If $w\in T$, then $e$ does not lie in $(G-e)[S]$. Hence $(G-e)[S]=G[S]$ has no perfect matching. By Lemma~\ref{lemma2026MinDegreeBasic}~\eqref{lemma2026MinDegreeBasic1}, $G[T]$ has no perfect matching. Hence $\overline{G}[T]$ is bad. Since $\overline{G-e}[T]\supseteq\overline{G}[T]$, the graph $\overline{G-e}[T]$ is bad, a contradiction. Thus $w\in S$, and similarly $x\in S$. So $e=wx$ and $S$ satisfy the hypotheses of Lemma~\ref{lemma2026MinDegreeMarkedge}, and at least one of \eqref{lemma2026MinDegreeMarkedge1}--\eqref{lemma2026MinDegreeMarkedge3} holds. If \eqref{lemma2026MinDegreeMarkedge1} holds, then $p,q$ are common neighbors of $w,x$ in $\overline{G}$, but $N_{\overline{G}}(w)\cap N_{\overline{G}}(x)=\emptyset$, a contradiction. If \eqref{lemma2026MinDegreeMarkedge2} holds, then $N_{\overline{G}[S]}(w)=N_{\overline{G}[S]}(x)\subseteq N_{\overline{G}}(w)\cap N_{\overline{G}}(x)=\emptyset$. Hence $d_{\overline{G}[S]}(w)=0$. But by Lemma~\ref{lemma2026MinDegreeMarkedge}, $d_{\overline{G}[S]}(w)\geq2$, a contradiction. Thus \eqref{lemma2026MinDegreeMarkedge3} holds.
    
            We have $\overline{G}[S]\cong K_{3,3}-wx$, with parts $\{w,a,b\}$ and $\{x,c,d\}$. Then $c,d\in N_{\overline{G}[S]}(w)\subseteq N_{\overline{G}}(w)=\{p,q\}$. So $\{c,d\}=\{p,q\}$ and the parts are $\{w,a,b\}$ and $\{x,p,q\}$. Since $\Delta(\overline{G})\leq3$ and the neighbors of $a$ in $\overline{G}[S]$ are $x,p,q$, we have $N_{\overline{G}}(a)=\{x,p,q\}$; similarly $N_{\overline{G}}(b)=\{x,p,q\}$. Hence every common neighbor of $a$ and $b$ is one of $x,p,q$. Since $wp,wq\in E(\overline{G})$ but $wx\notin E(\overline{G})$, the vertex $x$ is the only one of these three not adjacent to $w$ in $\overline{G}$. We now show that $w,p,q$ determine $x$. Since $p$ is adjacent to $w,a,b$ and has degree at most 3 in $\overline{G}$, $N_{\overline{G}}(p)=\{w,a,b\}$; similarly $N_{\overline{G}}(q)=\{w,a,b\}$. Hence $a$ and $b$ are fixed by $w,p,q$. Then $x$ is simply the third neighbor of $a$ besides $p,q$, and is also fixed. Therefore two different choices of $x$ would give the same $a,b$ and hence the same $x$, a contradiction. So among the five candidate vertices $x$, at most one can make $G-e$ not 3-extendable, and some choice makes $G-e$ 3-extendable.

        \item If $\delta(\overline{G})=3$, then $\overline{G}$ has 18 edges and $G$ has 48 edges. Call an edge $e=uv$ of $G$ \emph{obstructed} if $G-e$ is not 3-extendable. For an obstructed edge $e=uv$, by Lemma~\ref{lemma2026MinDegreeBasic}~\eqref{lemma2026MinDegreeBasic1}, there is a 6-set $T$ such that $(G-e)[T]$ has a perfect matching but $(G-e)[V(G)\setminus T]$ does not. Let $S=V(G)\setminus T$. Then $\overline{G-e}[S]$ is bad and $\overline{G-e}[T]$ is good. As in the case $\delta(\overline{G})=2$, we claim that $u,v\in S$. If $u\in T$, then $e$ does not lie in $(G-e)[S]$. So $(G-e)[S]=G[S]$ has no perfect matching. By Lemma~\ref{lemma2026MinDegreeBasic}~\eqref{lemma2026MinDegreeBasic1}, $G[T]$ has none. So $\overline{G}[T]$ is bad. Since $\overline{G-e}[T]\supseteq\overline{G}[T]$, the graph $\overline{G-e}[T]$ is also bad, a contradiction. Thus $e=uv$ and $S$ satisfy the hypotheses of Lemma~\ref{lemma2026MinDegreeMarkedge}, and one of \eqref{lemma2026MinDegreeMarkedge1}--\eqref{lemma2026MinDegreeMarkedge3} holds. Since $\delta(\overline{G})=3$ and $\Delta(\overline{G})\leq3$, then each vertex of $\overline{G}$ has degree exactly 3.
            
            If \eqref{lemma2026MinDegreeMarkedge1} holds, then $\overline{G}[\{u,v,p,q\}]\cong K_4-uv$. So $N_{\overline{G}}(p)=\{u,v,q\}$ and $N_{\overline{G}}(q)=\{u,v,p\}$. The edge $pq$ of $\overline{G}$ determines $e$. Distinct such edges correspond to distinct edges of $\overline{G}$. There are at most 18.
            
            If \eqref{lemma2026MinDegreeMarkedge2} holds, then $N_{\overline{G}}(u)=N_{\overline{G}}(v)$, both of size 3. Call two vertices \emph{twins} if they have the same neighborhood in $\overline{G}$. Then $u$ and $v$ are twins. Since a vertex does not belong to its own neighborhood, twins are never adjacent in $\overline{G}$. Hence every pair of twins are adjacent in $G$. We claim that each vertex has at most one twin. Suppose some vertex $u$ has two twins $v$ and $w$. So $N_{\overline{G}}(u)=N_{\overline{G}}(v)=N_{\overline{G}}(w)=\{a,b,c\}$. Then $u,v,w$ are pairwise non-adjacent in $\overline{G}$, and each of $a,b,c$ is adjacent to all of $u,v,w$. Since $\Delta(\overline{G})\leq3$, each of $a,b,c$ is adjacent to exactly $\{u,v,w\}$. In particular $a,b,c$ are pairwise non-adjacent. So $\overline{G}$ induces a $K_{3,3}$ on $\{u,v,w,a,b,c\}$ with parts $\{u,v,w\}$ and $\{a,b,c\}$, that is, a $K_{3,3}$-branch, a contradiction. Thus each vertex has at most one twin. Each such edge $uv$ is simply the edge of $G$ joining the twin pair $\{u,v\}$, and distinct twin pairs give distinct edges. Since each vertex has at most one twin, the 12 vertices split into at most 6 twin pairs. Hence there are at most 6 such edges.
            
            If \eqref{lemma2026MinDegreeMarkedge3} holds, then write the parts as $\{u,a,b\}$ and $\{v,c,d\}$. Since each vertex of $\overline{G}$ has degree exactly 3, $N_{\overline{G}}(a)=N_{\overline{G}}(b)=\{v,c,d\}$ and $N_{\overline{G}}(c)=N_{\overline{G}}(d)=\{u,a,b\}$. Thus $a$ and $b$ are twins, and so are $c$ and $d$. Now we count the obstructed edges (i.e. $uv$). In the structure of \eqref{lemma2026MinDegreeMarkedge3}, choose one endpoint of $uv$ from $\{v,c,d\}$. We have three choices. Once the endpoint, say $v$, is chosen, the other two common neighbors $c,d$ must be twins with common neighborhood $\{u,a,b\}$, where $u$ is another endpoint of $uv$ in the part containing $a,b$. So $u$ is the third neighbor of $c$ besides $a,b$. Hence for each choice there is at most one obstructed edge. So there are at most $3\cdot6=18$ obstructed edges by this way. Hence at most $18+6+18=42<48$ edges are obstructed. Some edge $e$ of $G$ is not obstructed, i.e., $G-e$ is 3-extendable.
    \end{itemize}

    In every case there is some edge $e$ such that $G-e$ is 3-extendable. So $G$ is not minimally 3-extendable, a contradiction.
\end{proof}\bigskip

\begin{proof}[Proof of Part \eqref{theorem2026MinDegreeK33} of Theorem~\ref{theorem2026MinDegreeK3}]
    By Lemma~\ref{lemma2026MinDegreePlummer}, $\delta(G)\geq4$. By Lemma~\ref{lemma2026MinDegree12Tight}, $\delta(G)\leq7$. If $\delta(G)=5$, then $\kappa(G)\leq5<6=2k$. By Lemma~\ref{lemma2026MinDegreeLouYu}~\eqref{lemma2026MinDegreeLouYu1}, $G$ is bipartite. Then by Lemma~\ref{lemma2026MinDegreeLouYu}~\eqref{lemma2026MinDegreeLouYu2}, we have $\delta(G)=4$. This is a contradiction. So $\delta(G)\in\{4,6,7\}$ and Part~\eqref{theorem2026MinDegreeK33} of Theorem~\ref{theorem2026MinDegreeK3} holds.
\end{proof}\bigskip

This completes the proof of Theorem~\ref{theorem2026MinDegreeK3}. Parts \eqref{theorem2026MinDegreeK31}, \eqref{theorem2026MinDegreeK32} and \eqref{theorem2026MinDegreeK33} are proved above. All predicted values occur. The value 4 occurs at order 8 (Proposition~\ref{proposition2026MinDegreeK44}) and at order 10 (Proposition~\ref{proposition2026MinDegreeK55}). The value 6 occurs at order 12 by $G_k$ with $k=3$ (Proposition~\ref{proposition2026MinDegreeMainconstr}). The value 7 occurs at order 8 (Proposition~\ref{proposition2026MinDegreeK8}) and at order 10 (Proposition~\ref{proposition2026MinDegreeTen}).

\section{Constructions}\label{constructions}

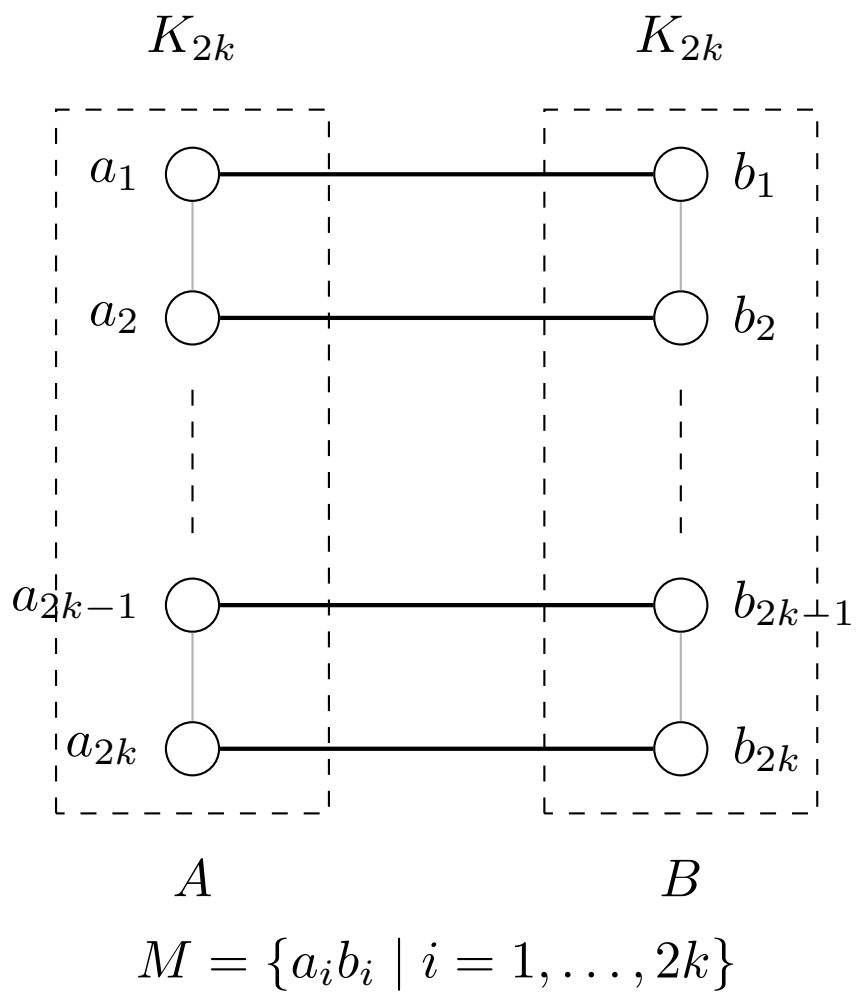
\begin{figure}[htbp]
    \centering
    \begin{tikzpicture}[scale=0.95,
      vtx/.style={circle, draw, fill=white, inner sep=0pt, minimum size=10pt, font=\small}]
      \node[vtx, label={[label distance=1pt]180:$a_1$}] (a1) at (0,2) {};
      \node[vtx, label={[label distance=1pt]180:$a_2$}] (a2) at (0,1) {};
      \node[vtx, label={[label distance=1pt]180:$a_{2k-1}$}] (aB1) at (0,-1) {};
      \node[vtx, label={[label distance=1pt]180:$a_{2k}$}] (aB2) at (0,-2) {};
      \node[vtx, label={[label distance=1pt]0:$b_1$}] (b1) at (3.4,2) {};
      \node[vtx, label={[label distance=1pt]0:$b_2$}] (b2) at (3.4,1) {};
      \node[vtx, label={[label distance=1pt]0:$b_{2k-1}$}] (bB1) at (3.4,-1) {};
      \node[vtx, label={[label distance=1pt]0:$b_{2k}$}] (bB2) at (3.4,-2) {};
    
      \draw[gray!60, thin] (a1) -- (a2);
      \draw[gray!60, thin] (aB1) -- (aB2);
      \draw[gray!60, thin] (b1) -- (b2);
      \draw[gray!60, thin] (bB1) -- (bB2);
    
      \draw[thick] (a1) -- (b1);
      \draw[thick] (a2) -- (b2);
      \draw[thick] (aB1) -- (bB1);
      \draw[thick] (aB2) -- (bB2);
     
      \draw[dashed] (0,0.5) -- (0,-0.5);
      \draw[dashed] (3.4,0.5) -- (3.4,-0.5);
      \draw[dashed] (-0.95,-2.45) rectangle (0.95,2.45);
      \draw[dashed] (2.45,-2.45) rectangle (4.35,2.45);
      \node at (0,2.95) {$K_{2k}$};
      \node at (3.4,2.95) {$K_{2k}$};
      \node at (0,-2.9) {$A$};
      \node at (3.4,-2.9) {$B$};
      \node at (1.7,-3.5) {$M=\{a_i b_i \mid i=1,\dots,2k\}$};
    \end{tikzpicture}
    \caption{A schematic of $G_k$. The \texttt{TikZ} code for this figure was prepared with the assistance of Hy3. The diagram itself was designed and verified by the authors.}
    \label{figureGk}
\end{figure}

We first give a construction that realizes the middle value $2k$ at the exceptional order $4k$, for every $k$. Let $k\geq1$. Take the vertex set $A\cup B$ with $A=\{a_1,\dots,a_{2k}\}$ and $B=\{b_1,\dots,b_{2k}\}$. Both $A$ and $B$ induce complete graphs. The remaining edges are $a_ib_i$ for $i=1,\dots,2k$. In other words, $G_k=\overline{K_{2k,2k}-M}$, where $M$ is a perfect matching of $K_{2k,2k}$. A schematic picture of $G_k$ is given in Figure~\ref{figureGk}.

\begin{proposition2026MinDegree}\label{proposition2026MinDegreeMainconstr}
    For every $k\geq1$ $G_k$ is a minimally $k$-extendable graph with $\nu(G_k)=4k$ and $\delta(G_k)=2k$.
\end{proposition2026MinDegree}
\begin{proof}
    Clearly we have $\nu(G_k)=4k$ and $\delta(G_k)=2k$.
    
    We show $G_k$ is $k$-extendable. Let $W$ be a matching of size $k$. Let $p,q,r$ be the numbers of its edges in $A$, in $B$, and of the form $a_ib_i$. Then $p+q+r=k$. An index $i$ is \emph{used} if $a_i$ or $b_i$ is an endpoint of an edge of $W$. The edges of $W$ in $A$ use $2p$ indices, those in $B$ use $2q$ indices, and the $r$ edges $a_ib_i$ use $r$ indices, so at most $2p+2q+r=2k-r$ indices are used and at least $r$ are unused. Exactly $2k-2p-r$ vertices of $A$ and $2k-2q-r$ vertices of $B$ are uncovered. Let $s=2k-2p-2q-r=r\geq0$. Choose $s$ unused indices $i$ and take the edges $a_ib_i$; these cover $s$ uncovered vertices on each side. The remaining $2q$ uncovered vertices of $A$ and $2p$ uncovered vertices of $B$ are paired into $q$ edges inside $A$ and $p$ edges inside $B$, which is possible since $A$ and $B$ are cliques. Together with $W$ this covers all $4k$ vertices. So $G_k$ is $k$-extendable.

    For minimality, the edges of $G_k$ may be represented by $a_1a_2$ and $a_1b_1$.

    Let $e=a_1a_2$ and $W=\{a_3a_4,\dots,a_{2k-1}a_{2k}\}\cup\{b_1b_2\}$, a matching of size $k$. The uncovered vertices are $a_1,a_2,b_3,\dots,b_{2k}$. In $G_k-e$ the vertex $a_1$ is adjacent to none of them. Its neighbors in $A$ are covered by $W$. Its only neighbor outside $A$ is $b_1$, covered by $W$. So $G_k-e$ is not $k$-extendable.

    Let $e=a_1b_1$ and $W=\{a_2b_2\}\cup\{a_3a_4,\dots,a_{2k-1}a_{2k}\}$, a matching of size $k$. The uncovered vertices are $a_1,b_1,b_3,\dots,b_{2k}$. In $G_k-e$ the vertex $a_1$ is adjacent to none of them. Its neighbors in $A$ are covered by $W$. So $G_k-e$ is not $k$-extendable. Hence $G_k$ is minimally $k$-extendable.
\end{proof}\bigskip

For $k=3$, this construction gives a minimally 3-extendable graph of order 12 and minimum degree 6. Our second construction gives a family of examples at order 10 with minimum degree 7.

\begin{proposition2026MinDegree}\label{proposition2026MinDegreeTen}
    Let $F$ be a triangle-free 2-regular graph on 10 vertices. Then $\overline{F}$ is minimally 3-extendable and $\delta(\overline{F})=7$.
\end{proposition2026MinDegree}
\begin{proof}
    We have $\delta(\overline{F})=7$. For every 4-set $S$, we have $\overline{F}[S]=K_4-F[S]$ and $F[S]$ is triangle-free with $\Delta\leq2$. By Lemma~\ref{lemma2026MinDegree4}, $\overline{F}[S]$ has a perfect matching. Hence, by Lemma~\ref{lemma2026MinDegreeBasic}~\eqref{lemma2026MinDegreeBasic2} with $k=3$ and $\nu=10$ (so $s=4$), $\overline{F}$ is 3-extendable. Let $e=xy$ be any edge of $\overline{F}$. Then $xy\notin E(F)$. Let $x',x''$ be the two neighbors of $x$ in $F$ and $S=\{x,y,x',x''\}$. Then $F[S]\cup\{e\}$ contains the star $K_{1,3}$ with center $x$. By Lemma~\ref{lemma2026MinDegree4}, the graph $(\overline{F}-e)[S]=K_4-(F\cup\{e\})[S]$ has no perfect matching. Note $V(\overline{F})\setminus S$ is a 6-set. Since $\Delta(F)\leq2$, by Corollary~\ref{lemma2026MinDegree6Cor}, $\overline{F}[V(\overline{F})\setminus S]=K_6-F[V(\overline{F})\setminus S]$ has a perfect matching. This is a 3-matching of $\overline{F}-e$ leaving exactly $S$, which has none. Hence $\overline{F}-e$ is not 3-extendable. So $\overline{F}$ is minimally 3-extendable.
\end{proof}\bigskip

The graphs $F$ are disjoint unions of cycles of length at least 4 summing to 10. They are $C_{10}$, $C_4\sqcup C_6$ and $C_5\sqcup C_5$. Finally, the value 4 also occurs at order 10.

\begin{proposition2026MinDegree}\label{proposition2026MinDegreeK55}
    Let $M$ be a perfect matching of $K_{5,5}$. Then $K_{5,5}-M$ is minimally 3-extendable and $\delta(K_{5,5}-M)=4$.
\end{proposition2026MinDegree}
\begin{proof}
    Let the bipartition be $A=\{a_1,\dots,a_5\}$, $B=\{b_1,\dots,b_5\}$ and $M=\{a_ib_i\}$. Then $a_ib_j$ is an edge exactly when $i\not=j$. A 3-matching uses three vertices from each part. It leaves $a_i,a_j$ and $b_p,b_q$ with $i\not=j$ and $p\not=q$. If $i\not=p$ and $j\not=q$, then choose the matching $\{a_ib_p,a_jb_q\}$. If $i=p$, then $i\not=q$ and $j\not= p$. Hence choose the matching $\{a_ib_q,a_jb_p\}$. Similarly for the case of $j=q$. So there is a matching covering the remaining vertices. So $G$ is 3-extendable. For any edge $e$, without loss of generality, we may assume $e=a_1b_2$. Let $W=\{a_2b_3,a_4b_5,a_5b_4\}$, a matching of size 3 disjoint from $e$. The uncovered vertices are $a_1,a_3,b_1,b_2$. Since in $G-e$ the vertex $a_1$ has no neighbor among these uncovered vertices, the only possible partners for $a_1$ are the non-edge $a_1b_1$ or the deleted edge $e=a_1b_2$. So $W$ cannot extend to any perfect matching in $G-e$. Hence $G-e$ is not 3-extendable.
\end{proof}

\section*{Statements and Declarations}

\noindent\textbf{Funding.}
The authors declare that no funds, grants, or other support were received during the preparation of this manuscript.

\noindent\textbf{Competing Interests.}
The authors have no relevant financial or non-financial interests to disclose.

\noindent\textbf{Author Contributions.}
Both authors contributed to the research. The first draft of the manuscript was written by Zhiyong Gan and both authors commented on previous versions of the manuscript. Both authors read and approved the final manuscript.

\noindent\textbf{Data Availability.}
Data sharing not applicable to this article as no datasets were generated or analysed during the current study.

\noindent\textbf{Use of AI Tools.}
During the preparation of this work the authors used Hy3 to assist in writing the \texttt{TikZ} code for Figure~\ref{figureGk}. The diagram was designed by the authors, and all its vertices and edges were checked against the definition of $G_k$. The authors reviewed and edited the content as needed and take full responsibility for the content of the published article.

\end{document}